# Bayesian model selection for linear regression

Miguel de Benito

Universität Augsburg

*Email:* `m.debenito.d@gmail.com`

Philipp Wacker

Universität Augsburg

*Email:* `phkwacker@gmail.com`

Abstract. In this note we introduce linear regression with basis functions in order to apply Bayesian model selection. The goal is to incorporate Occam's razor as provided by Bayes analysis in order to automatically pick the model optimally able to explain the data without overfitting.

Table of contents

## 1. Introduction: different models for linear regression

[1]Assume we have some collection of points $\{(x_n, t_n) \in \mathbb{R}^2 | n = 1, \ldots, N\}$. The $\boldsymbol{x} = (x_1, \ldots, x_N)$ are called **inputs** or **observations** and $\boldsymbol{t} = (t_1, \ldots, t_N)$ are the **outputs** or **targets**. One could

---

1. For a thorough treatment we refer to [Bis06, Chapter 3].

ask the following questions:

1. Is there any trend (monotonicity, asymptotics, etc.) in the data as a whole?

2. What's the likeliest value for the target $t$ given a new observation $x$?

3. Can we describe the target function as a linear combination of some set of known functions, at least approximately?

Here we will try to answer these questions using **linear regression**, which basically amounts to assuming the third item above: one first fixes some set of **basis functions** $\varphi_j \colon \mathbb{R} \to \mathbb{R}$, $j \in \{1, \ldots, M\}$ and then tries to fit scalar coefficients $w_1, \ldots, w_M \in \mathbb{R}$ such that for every $n \in \{1, \ldots, N\}$,

$$t_n \simeq \sum_{j=1}^{M} w_j \varphi_j(x_n)$$

in some suitable sense. A particular choice of basis functions is called an **hypothesis** or **model**, which we denote as $\mathcal{H}$ and it is an awesome feature of Bayesian statistics that one may let the data speak for themselves and tell us which one out of several hypothesis best fits the data, while at the same time not incurring in *overfitting*, i.e. the problem of choosing coefficients providing a perfect fit for a given training set but failing miserably with mostly any other one (this is done by skipping parameter fitting altogether, at least in principle). We will encode our lack of knowledge about what is the right model as a prior distribution on the set of all them (e.g. a uniform over the discrete set $\{\mathcal{H}_1, \ldots, \mathcal{H}_K\}$).

So, more generally than above, we start with several hypotheses $\mathcal{H}_1, \ldots, \mathcal{H}_K$ and try to choose the best one. To fix notation, we assume that the target function has the form

$$t = y(x, w) + \epsilon$$

with $\epsilon \sim \mathcal{N}(0, \sigma^2)$ a Gaussian random variable modeling noise, $x \in \mathbb{R}$ a data point and $w = (w_0, \ldots, w_{M-1}) \in W$ (e.g. $W = \mathbb{R}^M$), a vector of weights determining the function

$$y(x, w) = \sum_{j=0}^{M_k - 1} w_j \phi_j(x),$$

where

$$\phi_j = \phi_j^{(k)} \in \mathcal{H}_k = \{\phi_0^{(k)}, \ldots, \phi_{M_k - 1}^{(k)}\}, M_k \in \mathbb{N}$$

is some given set of basis functions. For convenience we always choose $\phi_0 \equiv 1$. One may use polynomials:

$$\mathcal{H}_1 = \{\phi_j(x) = x^j \colon j = 0, \ldots, M - 1\}, \quad M \in \mathbb{N},$$

or cosine functions

$$\mathcal{H}_2 = \{\phi_j(x) = \cos(\pi j x) \colon j = 0, \ldots, M - 1\}, \quad M \in \mathbb{N},$$

or wavelets or whatever, up to some fixed number $K$ of different models. The main task to perform now is **model selection**: given any pair of hypothesis $\mathcal{H}_{k_1}, \mathcal{H}_{k_2}$ determine which one is supported by more *evidence* in the data. In order to do this we will need to compute some integrals, and this we achieve using Laplace's method, which under our assumption of normality in the noise yields exact results. Finally, in order to apply this integration method we need to do **parameter inference**: for every hypothesis $\mathcal{H}_k$, compute (fit) the values for the weights which best explain the outputs $t$ in some suitable sense (which will be the maximization of the *a posteriori* distribution).

## 2. Parameter inference

**Question 1.** *Given an hypothesis $\mathcal{H}$ on how the output $t$ is generated from the observed data $x$ and some* a priori *distribution $p(w|\mathcal{H})$ on the set of parameters which govern the specific distribution, infer the correct posterior distribution $p(w|t, x, \mathcal{H})$ of the parameters given the data.*

**Notation.** *Since the observations $x$ are always given and we will not be doing any inference on their distribution or future values, we will omit the vector $x$ from the conditioning variables in the notation from here on. Thus, instead of $p(w|t, x, \mathcal{H})$ we will write $p(w|t, \mathcal{H})$, instead of $p(t|x, w, \mathcal{H})$ we will write $p(t|w, \mathcal{H})$ and so on. In a similar manner we will often refer to the outputs $t$ as the **data**, even though the whole data set includes $x$ as well.*

Solving this first question is a simple application of Bayes' rule:

$$p(w|t, \mathcal{H}) = \frac{p(t|w, \mathcal{H})}{p(t|\mathcal{H})} p(w|\mathcal{H}). \tag{1}$$

Because the normal is conjugate to itself, a convenient choice is $p(w|\mathcal{H}) = p(w|\text{no data}) \sim \mathcal{N}(\mu_w, \sigma_w^2 \text{ Id})$ and the new distribution will again be normal. We will be doing this later on, but other choices are possible! Notice that instead of just some numerical value for $w$, what we obtain using (1) is *the full a posteriori distribution*, which is a host more information, and thanks to the use of a conjugate prior we have an explicit expression.

As we have just seen, once we have all of the data $t$ we can compute the posterior in **batch** mode, i.e. using all of the data. However, it is remarkable that parameter inference can be done iteratively (**online**), i.e. updated when a new data point $t'$ arrives. The idea is that given the current dataset $t$, the corresponding posterior parameter distribution $p(w|t, \mathcal{H})$ is used as the *new prior* for $w$, when a new data point $t'$ is available. So we want to obtain the posterior distribution $p(w|t', t, \mathcal{H})$ as a function of $p(w|t, \mathcal{H})$, and this is essentially the same computation as above:

$$p(w|t', t, \mathcal{H}) = \frac{p(t'|t, w, \mathcal{H}) p(w|t, \mathcal{H})}{p(t'|t, \mathcal{H})} \stackrel{(*)}{=} \frac{p(t'|w, \mathcal{H})}{p(t'|\mathcal{H})} p(w|t, \mathcal{H}), \tag{2}$$

which again will be normal, and where (∗) is because of the independence of $t'$ wrt. $t$ given $w$. Also, $p(t)$ is constant in $w$ so we may ignore it for the computation of the law, as long as we normalize at the end.

Besides any possible reduction in computational cost, the real reason for doing parameter inference iteratively is that we can renormalize all terms after each step, so we will not run into machine precision problems. Consider equation (1): we will see in the next subsection that we do not need to calculate $p(t|\mathcal{H})$, so we only compute $p(t|w, \mathcal{H}) p(w|\mathcal{H})$, but this is essentially 0 for any realistic size of dataset $t$ and our parameter inference might fail. Done iteratively, we can divide by $p(t|\mathcal{H})$ at each step, which is just the normalization of the unnormalized density $p(t|w, \mathcal{H}) p(w|\mathcal{H})$ and all numbers should stay inside machine precision boundaries.

## 2.1. Method and the difference between ML and MAP.

From equations (1) or (2) we see that starting with a prior distribution $p(w|\mathcal{H})$ on the possible parameter values $w$ we can improve it using the acquired data $t$. It will be useful to define the following value for $w$. Let $W$ be the parameter space corresponding to hypothesis $\mathcal{H}$ (e.g. $W = \mathbb{R}^M$). The **maximum a posteriori** parameter $w_{\text{MAP}}$ is

$$w_{\text{MAP}} := \underset{w \in W}{\operatorname{argmax}}\, p(t|w, \mathcal{H}) p(w|\mathcal{H}) = \underset{w \in W}{\operatorname{argmax}}\, p(w|t, \mathcal{H}). \tag{3}$$

This is the *parameter value at which the posterior distribution of $w$ given the data attains its maximum.*[2] Bayesian reasoning uses the full distribution on $w$ given the data and the MAP value to account for the priors we choose (and priors that were posteriors of past updates, see figure 1).

$$
\begin{array}{ll}
p(w|\mathcal{H}) \longrightarrow & p(w|t, \mathcal{H}) = C \cdot p(t|w, \mathcal{H}) \cdot p(w|\mathcal{H}) \\
\text{initial prior} & \text{current posterior} = \text{(normalization)} \cdot \text{likelihood} \cdot \text{prior} \\
\\
p(w|t, \mathcal{H}) \longrightarrow & p(w|t, t', \mathcal{H}) = C' \cdot p(t'|w, \mathcal{H}) \cdot p(w|t, \mathcal{H}) \\
\text{new prior} & \text{current posterior} = \text{(normalization)} \cdot \text{likelihood} \cdot \text{prior}
\end{array}
$$

**Figure 1.** Bayesian updating: posteriors become new priors

It makes sense to take priors into account: if your data $t$ shows you consistent proof for a specific parameter value $w$, a single outlier $t'$ should not change your estimate much, i.e. we are not looking for the best fit but for the best fit consistent with our prior knowledge.

Finally note that we do not need to compute $p(t|\mathcal{H})$ since this is constant on $w$, as can be seen in the definition of $w_{\text{MAP}}$. More precisely, $p(t|\mathcal{H})$ is the normalization for the unnormalized density $w \mapsto p(t|w, \mathcal{H}) p(w|\mathcal{H})$.

---

2. This is not to be confused with the **maximum likelihood** parameter $w_{\text{ML}}$, which is defined as the *parameter value which best fits the data*: $w_{\text{ML}} := \operatorname{argmax}_{w \in W} p(t|w, \mathcal{H})$. Notice that $w_{\text{MAP}} \neq w_{\text{ML}}$ in general, especially when the best fit $w_{\text{ML}}$ is penalized by a small prior probability $p(w_{\text{ML}}|\mathcal{H})$ in the computation of $w_{\text{MAP}}$.

## 3. Model selection

In order to compare two models $\mathcal{H}_1, \mathcal{H}_2$ given a dataset $t = \{t_1, \ldots, t_N\}$, we want to compute the ratio of the hypotheses' probability, i.e.

$$\frac{p(\mathcal{H}_1|t)}{p(\mathcal{H}_2|t)} = \frac{p(t|\mathcal{H}_1)}{p(t|\mathcal{H}_2)} \frac{p(\mathcal{H}_1)}{p(\mathcal{H}_2)}. \tag{4}$$

Notice that, as in the case of batch parameter estimation this introduces a problem of machine precision. For any decent number of data points, any given dataset will have almost vanishing probability. This is true of any probability measure we choose but is most obviously seen e.g. if the data comes from the realization of i.i.d. random variables, since then $p(t|\mathcal{H}_k) = \prod_{n=1}^{N} p(t_n|\mathcal{H}_k) \ll 1$. Therefore the quotient above is roughly 0/0 and we have a problem. In parameter inference we solved that problem by doing iterative steps and renormalizing each time, but as we will see, this can not be done for model selection in general, at least in a naive way:

### 3.1. Why (naive) iterative model selection does not work.

The naive idea one (for example the authors of this manuscript) could have would be to implement an online version of model selection mimicking what we did in parameter estimation. Given a dataset $t = (t_1, \ldots, t_N)$ and a new data point $t' = t_{N+1}$ we compute

$$p(\mathcal{H}_n|t, t') = \frac{p(t'|\mathcal{H}_k, t) p(\mathcal{H}_k|t)}{p(t'|t)} \stackrel{(*)}{=} \frac{p(t'|\mathcal{H}_k)}{p(t')} p(\mathcal{H}_k|t), \tag{5}$$

where $(*)$ would be deemed correct due to the independence of $t'$ wrt. $t$, as, by independence:

$$p(t'|\mathcal{H}_k, t) = \frac{p(t', t|\mathcal{H}_k)}{p(t|\mathcal{H}_k)} = \frac{p(t'|\mathcal{H}_k) p(t|\mathcal{H}_k)}{p(t|\mathcal{H}_k)} = p(t'|\mathcal{H}_k). \tag{6}$$

**Warning.** Formulas (5) and (6) are wrong! For a intuitive example on why this is so, see section 4.3.

The root of all evil is that

$$p(t, t'|r, \mathcal{H}_k) = p(t|r, \mathcal{H}_k) p(t'|r, \mathcal{H}_k)$$

but

$$p(t, t'|\mathcal{H}_k) \neq p(t|\mathcal{H}_k) p(t'|\mathcal{H}_k),$$

i.e. independence only holds conditional to a specific parameter value $r$.

Put in another way: when we decided that we would have a prior on the parame-ters and that the distribution of the observations would be given by conditional probabilities, we were saying that the joint distribution is

$$p(\boldsymbol{t}, w) = p(w) \prod_{n=1}^{N} p(t_n|w). \qquad (7)$$

This setting can be depicted as in Figure 2.

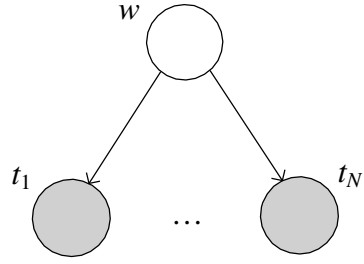

**Figure 2.** Bayesian net for i.i.d. data

Using now the product rule in (7), we have $p(\boldsymbol{t}|w) = \prod_{n=1}^{N} p(t_n|w)$, that is: the observations are **conditionally independent** given the parameters (one writes $t_n \perp t_m | w$). However, if dispense with $w$ by marginalizing $p(\boldsymbol{t}, w) = p(\boldsymbol{t}|w) p(w)$ over it we obtain

$$p(\boldsymbol{t}) = \int_W p(\boldsymbol{t}|w) p(w) \, dw \neq \prod_{n=1}^{N} p(t_n).$$

### 3.2. Correct model selection.

The proper way of doing model selection for a new data point involves considering the whole data $(\boldsymbol{t}, t')$ instead of using just the new datum $t'$ to update the prior. For a feeling of why this is necessary see the remarks in section 4.3. As we can't do model selection iteratively, we will state the model selection formula as depending on the whole data $\boldsymbol{t}$:

$$p(\mathcal{H}_k | \boldsymbol{t}) = \frac{p(\boldsymbol{t}|\mathcal{H}_k) p(\mathcal{H}_k)}{p(\boldsymbol{t})} = \frac{\int_W p(\boldsymbol{t}|w, \mathcal{H}_k) p(w|\mathcal{H}_k) \, dw}{p(\boldsymbol{t})} p(\mathcal{H}_k). \qquad (8)$$

We can see in equation (8) that the a posteriori probability for $\mathcal{H}_k$ only "improves" if the data $\boldsymbol{t}$ are better explained by the specific hypothesis $\mathcal{H}_k$ than they would be by just using the prior, i.e. if the ratio $p(\boldsymbol{t}|\mathcal{H}_k)/p(\boldsymbol{t})$ is greater than one. The quantities in equation (8) are:

- $p(\mathcal{H}_k)$ is the prior distribution for the $\mathcal{H}_k$ (e.g. uniform).

- $p(\boldsymbol{t})$ is a quantity we won't have to compute (see below).

- $p(\boldsymbol{t}|\mathcal{H}_k)$ is the probability of the datum $\boldsymbol{t}$ being explained by *any possible choice of parameters in the model $\mathcal{H}_k$*. We call it the **evidence for $\mathcal{H}_k$**. We must integrate over parameter space in order to compute it:

$$p(\boldsymbol{t}|\mathcal{H}_k) = \int_W p(\boldsymbol{t}|w, \mathcal{H}_k) p(w|\mathcal{H}_k) \, dw. \qquad (9)$$

For the same reasons as in parameter inference, the denominator (in this case $p(\boldsymbol{t})$) can be ignored: it is a normalization factor constant between all hypotheses. Hence we only need to find the numerator $p(\boldsymbol{t}|\mathcal{H}_k) p(\mathcal{H}_k) = \int p(\boldsymbol{t}|w, \mathcal{H}_k) p(w|\mathcal{H}_k) \, dw \cdot p(\mathcal{H}_k)$. Then we can normalize by $p(\boldsymbol{t}) = \sum_k p(\boldsymbol{t}|\mathcal{H}_k) p(\mathcal{H}_k)$.

**Remark 2.** The integral (9) is in general hard to compute. We could try to use MCMC, but given our assumptions of normality it is best to use Laplace's Method discussed below.

### 3.3. Laplace's method.

There is a way of tackling the integral in equation (9) and the denominator in the update term of parameter estimation in equation (2). You can read the details in [Mac05] (especially the great section about the interpretation of the "Occam factor"!), we will only collect the results:

$$p(t|\mathcal{H}_k) = \int_W p(t|w, \mathcal{H}_k) p(w|\mathcal{H}_k) \, dw \approx \frac{p(t|w_{\text{MAP}}, \mathcal{H}_k) p(w_{\text{MAP}}|\mathcal{H}_k)}{\sqrt{\det(A/2\pi)}}, \quad (10)$$

where $w_{\text{MAP}} = w_{\text{MAP}}^k$ is the maximum a posteriori parameter defined in (3) and

$$A := -\nabla_w^2 \ln p(w|t, \mathcal{H}_k)\big|_{w=w_{\text{MAP}}}$$
$$\stackrel{(2)}{=} [-\nabla_w^2 \ln p(t|w, \mathcal{H}_k) - \nabla_w^2 \ln p(w|\mathcal{H}_k)]\big|_{w=w_{\text{MAP}}}. \quad (11)$$

## 4. An easy example: is this coin bent?

Consider a coin whose fairness is in doubt. We would like to infer the probability of it being fraudulent as well as its "bentness".

### 4.1. Model.

We will call the sides of the coin 1 and 0. The result of flipping the coin is a random event independent from any past events. Our rivaling hypotheses are

- $\mathcal{H}_0$: The coin is fair, i.e. $\mathbb{P}(0) = \mathbb{P}(1) = 1/2$.

- $\mathcal{H}_1$: The coin has been manipulated: $\mathbb{P}(0) = r$ and $\mathbb{P}(1) = 1 - r$. We assume an initial non-committal distribution on $r$, i.e. the density of $r$ is $p(r|\mathcal{H}_1) = 1$ on $[0,1]$.

Initially, we consider both hypotheses as being identically probable, i.e. $\mathbb{P}(\mathcal{H}_k) = 1/2$.

**Remark 3.** Let us interpret the probabilities as statistics: Our prior on the hypotheses tells us that the experiment's conductor thinks that it is equally likely that he gets a bent coin or a fair coin. A person unfamiliar with casinos might use a different prior, e.g. $\mathbb{P}(\mathcal{H}_0) = 0.99$. If your beliefs are very strong at the beginning, the amount of data you need to change your mind is very large.

The prior on the parameter $r$ is even more strongly debatable: assuming that a coin has been manipulated, we could think that the probability of $r$ being 0.5 should be 0. Also, if the manipulator does not want his forgery to be too easily detected, he is likely to be more prudent than to set $r = 0$ or $r = 1$. Perhaps a double-peak distribution with maxima in $r = 0.5 \pm 0.03$ might be a better prior (just out of experience and reasoning). For easier calculation, we will use a uniform prior on $r$ but beware that Bayesian inference is all about your priors and beliefs about reality. Bayesian statistics states that *there is no such thing as objective inference*, just a correct propagation of belief states. In practice, in a lot of cases, different priors only change the speed of convergence of your inference and unless your initial distribution deems ranges of events as impossible (i.e. $p(r|\mathcal{H}_1)=0$ for some $r$), collecting more and more data will yield similar conclusions between different priors.

Also note that $\mathcal{H}_1$ is a perfect super-hypothesis of $\mathcal{H}_0$ in the sense that $\mathcal{H}_0$ is included in $\mathcal{H}_1$ (for $r=0.5$).

As a function of the model (and its parameter, if it has any), we get the following coin event probability distributions:

- $p(t=0|\mathcal{H}_0)=\frac{1}{2}=p(t=1|\mathcal{H}_0)$
- $p(0|r,\mathcal{H}_1)=r=:r_0$ and $p(1|r,\mathcal{H}_1)=1-r=:r_1$

### 4.2. Parameter inference.

Now we try to infer the parameter $r$ in model $\mathcal{H}_1$ by throwing the coin and drawing conclusions from the outcome. Assume that we are in the middle of inference and we have the distribution of $r$ given the data $\boldsymbol{t} = (t_1, \ldots, t_N)$. As soon as the next data point $t'=t_{N+1}\in\{0,1\}$ arrives, we have:

$$\begin{aligned}
p(r|\boldsymbol{t},t',\mathcal{H}_1) &= \frac{p(t'|r,\mathcal{H}_1)\,p(r|\boldsymbol{t},\mathcal{H}_1)}{p(t'|\mathcal{H}_1)} \\
&= \frac{r_{t'}}{\int_0^1 p(t'|\rho,\mathcal{H}_1)\,p(\rho|\mathcal{H}_1)\,\mathrm{d}\rho}\,p(r|\boldsymbol{t},\mathcal{H}_1) \\
&= \begin{cases}
\dfrac{r}{\int_0^1 \rho\,\mathrm{d}\rho}\,p(r|\boldsymbol{t},\mathcal{H}_1) \;=\; 2r\,p(r|\boldsymbol{t},\mathcal{H}_1) & \text{if } t'=0, \\[1ex]
\dfrac{1-r}{\int_0^1 (1-\rho)\,\mathrm{d}\rho}\,p(r|\boldsymbol{t},\mathcal{H}_1) \;=\; 2(1-r)\,p(r|\boldsymbol{t},\mathcal{H}_1) & \text{if } t'=1.
\end{cases}
\end{aligned}$$

For future use we note the following:

$$p(t'|\mathcal{H}_1)=1/2=p(t'|\mathcal{H}_0). \tag{12}$$

Now that we have the updated parameter distribution we are interested in which model might be the most probable given the data.

### 4.3. Iterative model selection. Why it cannot work.

If we try to infer the model iteratively from the data according to the (wrong) formula (5), we obtain using (12):

$$\frac{p(\mathcal{H}_0|\boldsymbol{t},t')}{p(\mathcal{H}_1|\boldsymbol{t},t')} = \frac{p(t'|\mathcal{H}_0)}{p(t'|\mathcal{H}_1)}\frac{p(\mathcal{H}_0|\boldsymbol{t})}{p(\mathcal{H}_1|\boldsymbol{t})} \stackrel{(12)}{=} \frac{p(\mathcal{H}_0|\boldsymbol{t})}{p(\mathcal{H}_1|\boldsymbol{t})} = \ldots = \frac{p(\mathcal{H}_0)}{p(\mathcal{H}_1)},$$

which means that the result of model selection yields a prior independent from the data and this is clearly wrong. As we saw in Section 3.1, the problem is that $p(t, t'|\mathcal{H}_1) \neq p(t|\mathcal{H}_1) p(t'|\mathcal{H}_1)$. To illustrate this consider $t = (1, 1, 1, \ldots, 1)$ (say 100 times). Then $r$ will very probably be approximately 0 and $p(t, t'|\mathcal{H}_1) \approx p(t, t'|r=0, \mathcal{H}_1) = 1\, p(t|\mathcal{H}_1)$, but $p(t|\mathcal{H}_1) p(t'|\mathcal{H}_1) = 0.5\, p(t|\mathcal{H}_1)$. The idea is that the occurrence of $t$ might set a strong bias for the possible values of $r$ which in turn change the probability distribution of $t'$.

### 4.4. Model selection: "batch" approach.

The way to select the correct model is

$$p(\mathcal{H}_1|t) = \frac{p(t|\mathcal{H}_1) p(\mathcal{H}_1)}{p(t)} = \frac{\int p(t|r, \mathcal{H}_1) p(r|\mathcal{H}_1) \, dr}{p(t)} p(\mathcal{H}_1)$$

and

$$p(\mathcal{H}_0|t) = \frac{p(t|\mathcal{H}_0)}{p(t)} p(\mathcal{H}_0).$$

The key is getting the evidences $p(t|\mathcal{H}_1) = \int p(t|r, \mathcal{H}_1) p(r|\mathcal{H}_1) \, dr$ and $p(t|\mathcal{H}_0)$, the remaining factors are the model priors (which are known) and $p(t)$ (which is just the normalization). The evidence for $\mathcal{H}_0$ is easy: consider $t = (t_1, \ldots, t_N)$, then

$$p(t|\mathcal{H}_0) = 2^{-N}.$$

To compute the evidence for $\mathcal{H}_1$ define $h(t) = $ ``amount of 0s in $t$''. Now $p(t|r, \mathcal{H}_1) = r^{h(t)} (1-r)^{N-h(t)}$ and using

$$\int_0^1 r^k (1-r)^{N-k} \, dr = \frac{k!(N-k)!}{(N+1)!} = \frac{1}{(N+1)\binom{N}{k}},$$

we get

$$\frac{p(\mathcal{H}_1|t)}{p(\mathcal{H}_0|t)} = \frac{\int_0^1 r^{h(t)}(1-r)^{N-h(t)} \, dr}{2^{-N}} \frac{p(\mathcal{H}_1)}{p(\mathcal{H}_0)} = \frac{p(\mathcal{H}_1) 2^N}{p(\mathcal{H}_0)(N+1)\binom{N}{h(t)}} \tag{13}$$

This shows a classic problem of this kind of methods: the numbers involved grow quickly and exact computation becomes unfeasible.

Consider the example of a fair coin: For the case $h(t) \simeq N/2$ use Stirling's formula $n! \simeq \sqrt{2\pi n}\, (n/e)^n$ and see

$$\frac{2^N}{(N+1)\binom{N}{h(t)}} \simeq \frac{2^N \left[\sqrt{2\pi \frac{N}{2}} \left(\frac{N}{2e}\right)^{\frac{N}{2}}\right]^2}{\sqrt{2\pi \frac{N+1}{2}} \left(\frac{N+1}{e}\right)^{N+1}}$$

$$\simeq \frac{\sqrt{\pi}}{\sqrt{N+1}}$$

$$\simeq O\!\left(\frac{1}{\sqrt{N}}\right),$$

so the hypothesis $\mathcal{H}_0$ becomes probable with speed $N^{-1/2}$ if there is no specific evidence for $\mathcal{H}_1$.

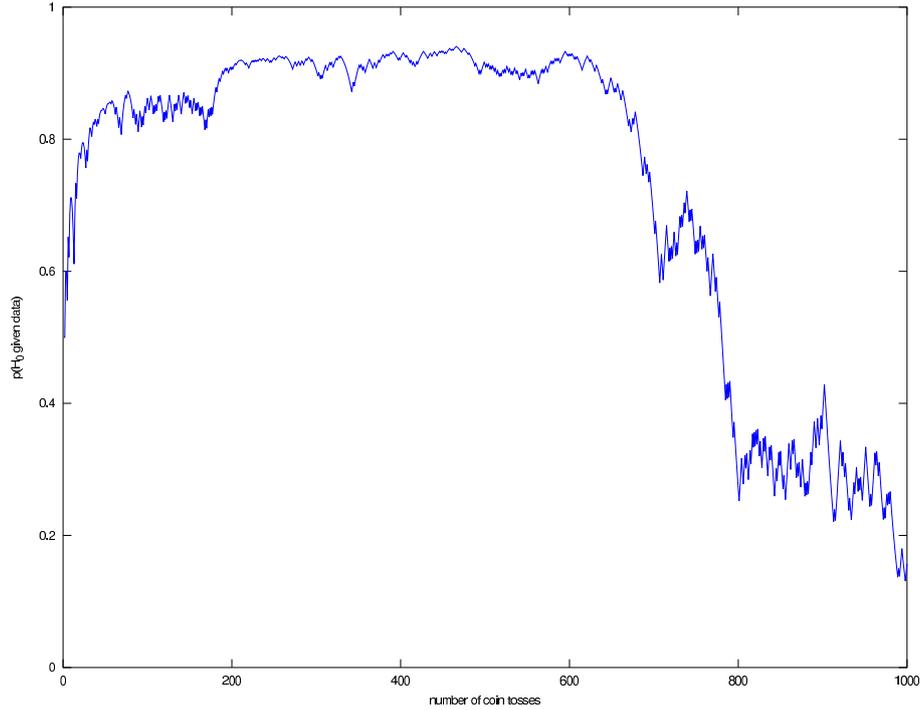

**Figure 3.** Model selection for a bent coin ($p(\text{``Heads"})=0.55$) for 1000 tosses. Initially, the slight imbalance between the number of heads and tails can be explained perfectly by a fair coin, so this hypothesis is strongly supported. After 600 tosses, the probability that the coin is in fact bent gets more and more support.

### 4.5. Model selection: Quasi-iterative approach.

We showed in sections 3.1 and 4.3 why iterative methods cannot work in a strict sense. There is a way which still uses only computationally feasible terms by taking a clever updating approach possible because one can explicitly compute the integrals involved. Let $t = (0,0,1,0,1,0,0,\ldots)$ be a data vector with 0 being heads and 1 being tails. Assume that the length of $t$ is $N$ and the number of 0s is $K$. Now a new data point $t'$ arrives, which is either 0 or 1. We write $K'$ for the number of 0s after arrival of $t'$, so $K' \in \{K, K+1\}$. Then

$$\frac{p(\mathcal{H}_0|t,t')}{p(\mathcal{H}_1|t,t')} = \frac{p(t'|\mathcal{H}_0,t)}{p(t'|\mathcal{H}_1,t)} \cdot \frac{p(\mathcal{H}_0|t)}{p(\mathcal{H}_1|t)} \tag{14}$$

The last term is the prior version of the left hand side fraction. Thus, the term $\frac{p(t'|\mathcal{H}_0,t)}{p(t'|\mathcal{H}_1,t)}$ is an update term we can use. First,

$$p(t'|\mathcal{H}_0,t) = p(t'|\mathcal{H}_0) = \frac{1}{2}.$$

Then

$$p(t'|\mathcal{H}_1,t) = \frac{p(t',t|\mathcal{H}_1)}{p(t|\mathcal{H}_1)} \tag{15}$$

We expand to calculate

$$p(t|\mathcal{H}_1) = \int_0^1 p(t|r,\mathcal{H}_1)\,p(r|\mathcal{H}_1)\,dr = \int_0^1 r^K(1-r)^{N-K}\,dr = \frac{1}{(N+1)\binom{N}{K}}$$

and

$$p(t',t|\mathcal{H}_1) = \int_0^1 r^{K'}(1-r)^{N+1-K'}\,dr = \frac{1}{(N+2)\binom{N+1}{K'}}.$$

Combining this with (15) and (14), we obtain

$$\frac{p(\mathcal{H}_0|t,t')}{p(\mathcal{H}_1|t,t')} = \begin{cases} \dfrac{N+2}{2K+2} \dfrac{p(\mathcal{H}_0|t)}{p(\mathcal{H}_1|t)} & \text{if } t'=0, \\ \dfrac{N+2}{2(N+1-K)} \dfrac{p(\mathcal{H}_0|t)}{p(\mathcal{H}_1|t)} & \text{if } t'=1. \end{cases}$$

This means that we do not need to save the exact data vector but just the number of heads and tails to apply a recursive-type model selection approach. Also, in comparison with formula (13), we will not run into computational overflow as there are no exponentially growing terms.

## 5. Influence of the race on the imposition of the death penalty

Consider exercise 28.4 in [Mac05]. Here Occam's razor yields a tool of analyzing whether there is evidence for racial discrimination in the data given. Consider the events

- $V$: The victim is white (and $V^C$: The victim is black)
- $M$: The defendant is white (and $M^C$: The defendant is black)
- $D$: The defendant gets sentenced to death

The hypotheses we will look at are

- $\mathcal{H}_{00}$: Neither the victim's nor the defendant's race has influence on the probability of the defendant getting the death penalty, i.e.
  $$p(D) = \tau.$$
- $\mathcal{H}_{10}$: Only the victim's race influences the death penalty probability:
  $$p(D|V) = \tau, \quad p(D|V^C) = \chi.$$
- $\mathcal{H}_{01}$: Only the defendant's race influences the death penalty probability:
  $$p(D|M) = \tau, \quad p(D|M^C) = \chi.$$
- $\mathcal{H}_{11}$: Both victim's and defendant's race influence the death penalty probability:
  $$p(D|V,M) = \tau, \quad p(D|V^C,M) = \chi, \quad p(D|V,M^C) = \rho, \quad p(D|V^C,M^C) = \sigma.$$

We will always assume uniform priors of all parameters $\tau$ etc. for every hypothesis.

For convenience we start with the log probabilities:

$$\begin{aligned}
\ln p(t|\tau, \mathcal{H}_{00}) &= (19+11+6)\ln(\tau) + (132+9+52+97)\ln(1-\tau) \\
\ln p(t|\tau, \chi, \mathcal{H}_{10}) &= (19+11)\ln(\tau) + (132+52)\ln(1-\tau) \\
&\quad + 6\ln(\chi) + (9+97)\ln(1-\chi) \\
\ln p(t|\tau, \chi, \mathcal{H}_{01}) &= 19\ln(\tau) + (132+9)\ln(1-\tau) \\
&\quad + (11+6)\ln(\chi) + (52+97)\ln(1-\chi) \\
\ln p(t|\tau, \chi, \rho, \sigma, \mathcal{H}_{11}) &= 19\ln(\tau) + 132\ln(1-\tau) \\
&\quad + 0\ln(\chi) + 9\ln(1-\chi) \\
&\quad + 11\ln(\rho) + 52\ln(1-\rho) \\
&\quad + 6\ln(\sigma) + 97\ln(1-\sigma).
\end{aligned}$$

As our priors are uniform, the MAP parameters are the ML parameters which in turn are

- $\mathcal{H}_{00}$: $\tau_{\text{MAP}} = \frac{19+11+6}{19+11+6+132+9+52+97} \approx 0.1104$.

- $\mathcal{H}_{10}$: $\tau_{\text{MAP}} = \frac{19+11}{19+11+132+52} \approx 0.1402$ and $\chi_{\text{MAP}} = \frac{6}{6+9+97} \approx 0.0536$.

- $\mathcal{H}_{01}$: $\tau_{\text{MAP}} = \frac{19}{19+132+9} \approx 0.1187$ and $\chi_{\text{MAP}} = \frac{11+6}{11+6+52+97} \approx 0.1024$.

- $\mathcal{H}_{11}$: $\tau_{\text{MAP}} = \frac{19}{19+132} \approx 0.1258$, $\chi_{\text{MAP}} = 0$, $\rho_{\text{MAP}} = \frac{11}{11+52} \approx 0.1746$ and $\sigma_{\text{MAP}} = \frac{6}{6+97} \approx 0.0583$.

One can observe that when the victim was white and the murderer black, there was a 17.46 percent probability of the defendant being sentenced to death whereas the opposite combination did not lead to any single death penalty sentence.

Recall ($w$ being the abstract parameter variable)

$$p(t|\mathcal{H}_k) = \int p(t|w,\mathcal{H}_k) p(w|\mathcal{H}_k) dw \approx \frac{p(t|w_{\text{MAP}}, \mathcal{H}_k) p(w_{\text{MAP}}|\mathcal{H}_k)}{\sqrt{\det(A/2\pi)}}$$

where $A = -\nabla_w^2 \ln p(t|w,\mathcal{H}_k) - \nabla_w^2 \ln p(w|\mathcal{H}_k)$, where the latter term is zero in our case as we have uniform priors. The individual matrices $A$ are:

$$A_{00} = \frac{36}{\tau^2} + \frac{290}{(1-\tau)^2}$$

$$A_{10} = \text{diag}\left(\frac{30}{\tau^2} + \frac{184}{(1-\tau)^2}, \frac{6}{\chi^2} + \frac{106}{(1-\chi)^2}\right)$$

$$A_{01} = \text{diag}\left(\frac{19}{\tau^2} + \frac{141}{(1-\tau)^2}, \frac{17}{\chi^2} + \frac{149}{(1-\chi)^2}\right)$$

$$A_{11} = \text{diag}\left(\frac{19}{\tau^2} + \frac{132}{(1-\tau^2)}, \frac{9}{(1-\chi)^2}, \frac{11}{\rho^2} + \frac{52}{(1-\rho)^2}, \frac{6}{\sigma^2} + \frac{97}{(1-\sigma)^2}\right).$$

Then the evidence for each model is (read subscript "MAP" for each parameter)

$$p(t|\mathcal{H}_{00}) \approx \frac{\tau^{36}(1-\tau)^{290}}{\sqrt{\det(A_{00}/2\pi)}}$$
$$\approx 2.8313 \cdot 10^{-51}$$
$$p(t|\mathcal{H}_{10}) \approx \frac{\tau^{30}(1-\tau)^{184} \chi^6 (1-\chi)^{106}}{\sqrt{\det(A_{10}/2\pi)}}$$
$$\approx 4.698 \cdot 10^{-51}$$
$$p(t|\mathcal{H}_{01}) \approx \frac{\tau^{19}(1-\tau)^{141} \chi^{17} (1-\chi)^{149}}{\sqrt{\det(A_{01}/2\pi)}}$$
$$\approx 2.7485 \cdot 10^{-51}$$
$$p(t|\mathcal{H}_{11}) \approx \frac{\tau^{19}(1-\tau)^{132} \chi^0 (1-\chi)^9 \rho^{11}(1-\rho)^{52} \sigma^6 (1-\sigma)^{97}}{\sqrt{\det(A_{11}/2\pi)}}$$
$$\approx 1.4875 \cdot 10^{-51}.$$

As $\sum_{i,j=0}^{1} p(\mathcal{H}_{i,j}|t) = 1$, and assuming non-committal priors $p(\mathcal{H}_{i,j}) = \frac{1}{4}$ for $i, j = 0, 1$, we obtain

$$p(\mathcal{H}_{00}|t) \approx 0.24,$$
$$p(\mathcal{H}_{10}|t) \approx 0.40,$$
$$p(\mathcal{H}_{01}|t) \approx 0.23,$$
$$p(\mathcal{H}_{11}|t) \approx 0.13.$$

Hence, the hypothesis $\mathcal{H}_{10}$ stating that the race of the victim is contributing to the probability of a murderer getting sentenced to death is the most probable from a Bayesian point of view *given the priors we had.*

## 6. Bayesian Model Selection for Linear Regression

### 6.1. Setting.

We will consider the following polynomial hypotheses:

$$\mathcal{H}_k^P: \quad t_n = w_0 + w_1 x_n + \cdots + w_{k-1} x_n^{k-1} + \varepsilon_n, \quad n = 1, \ldots, N, k = 1, \ldots, K.$$

For the following we can just as well use different hypothesis spaces as for example trigonometric functions:

$$\mathcal{H}_k^T: \quad t_n = \sum_{j=1}^{k} \cos(j \pi x_n) + \sin(j \pi x_n) + \varepsilon_n, \quad n = 1, \ldots, N, k = 1, \ldots, K$$

For the sake of concreteness, we will state everything for the example of polynomial functions, but everything can be generalized to generic function spaces.

Let $\varepsilon_n \sim \mathcal{N}(0, \sigma^2)$ i.i.d. where $\sigma$ is a known noise parameter (for this application we will not infer the noise from the data although this can be done). For brevity we denote

$$\Phi_k(x_n) = (1, x, x^2, \ldots, x^{k-1})$$

and the full data matrix is

$$\Phi_k = \begin{pmatrix} \Phi_k(x_1) \\ \vdots \\ \Phi_k(x_N) \end{pmatrix} = \begin{pmatrix} 1 & x_1 & x_1^2 & \cdots & x_1^{k-1} \\ 1 & x_2 & x_2^2 & \cdots & x_2^{k-1} \\ \vdots & \vdots & \vdots & \ddots & \vdots \\ 1 & x_N & x_N^2 & \cdots & x_N^{k-1} \end{pmatrix}.$$

Then

$$p(t_n | \boldsymbol{w}, \mathcal{H}_k) = \mathcal{N}(\boldsymbol{w}^\top \Phi_k(x_n)^\top, \sigma^2).$$

We assume the prior model probabilities to be equal, i.e. $p(\mathcal{H}_k) = 1/K$. The parameters will be distributed normally, i.e. the density of the parameter vector $\boldsymbol{w}$ is

$$p(\boldsymbol{w} | \mathcal{H}_k) = \mathcal{N}(\boldsymbol{0}, \sigma_w^2 \mathrm{Id}_k).$$

### 6.2. Iterative parameter estimation.

As explained in Section 2, to avoid numerical problems we will proceed iteratively. In step $n > 1$ we use the parameter estimate given by the posterior $p(\boldsymbol{w}|t_1, \ldots, t_{n-1}, \mathcal{H}_k)$ coming from step $n-1$ as the new prior, and for $n = 1$ we have the prior $p(\boldsymbol{w}|\mathcal{H}_k)$. The likelihood for the new data point $t_n$ is also known:

$$p(t_n|\boldsymbol{w}, \mathcal{H}_k) = \mathcal{N}(\boldsymbol{w}^\top \Phi_k(x_n), \sigma^2).$$

Then

$$p(\boldsymbol{w}|t_n, \mathcal{H}_k) = \frac{p(t_n|\boldsymbol{w}, \mathcal{H}_k)\, p(\boldsymbol{w}|\mathcal{H}_k)}{p(t_n|\mathcal{H}_k)}$$

and (cf. (3))

$$\boldsymbol{w}_{\mathrm{MAP}} = \underset{\boldsymbol{w} \in W_k}{\operatorname{argmax}}\, p(\boldsymbol{w}|t_n, \mathcal{H}_k) = \underset{\boldsymbol{w} \in W_k}{\operatorname{argmax}}\, p(t_n|\boldsymbol{w}, \mathcal{H}_k)\, p(\boldsymbol{w}|\mathcal{H}_k).$$

As all distributions are Gaussian, the parameter estimation is easy (see for example [Bis06]) and will not be done here.

### 6.3. Model selection.

As explained in Section 3, model selection works by ordering the hypotheses according to their a posteriori probability

$$p(\mathcal{H}_k|\boldsymbol{t}) = \frac{p(\boldsymbol{t}|\mathcal{H}_k)}{p(\boldsymbol{t})}\, p(\mathcal{H}_k).$$

This is done by ranking their evidence $p(\boldsymbol{t}|\mathcal{H}_k)$, computed using Laplace's method (see Section 3.3) as

$$p(\boldsymbol{t}|\mathcal{H}_k) = \int_{W_k} p(\boldsymbol{t}|\boldsymbol{w}, \mathcal{H}_k)\, p(\boldsymbol{w}|\mathcal{H}_k)\, d\boldsymbol{w} \approx \frac{p(\boldsymbol{t}|\boldsymbol{w}_{\mathrm{MAP}}, \mathcal{H}_k)\, p(\boldsymbol{w}_{\mathrm{MAP}}|\mathcal{H}_k)}{\sqrt{\det(A/2\pi)}},$$

where $\boldsymbol{w}_{\mathrm{MAP}}$ is given by (3) and $A$ is given by (11). The log-likelihood of the data is

$$\ln p(\boldsymbol{t}|\boldsymbol{w}, \mathcal{H}_k) = C - \frac{\sum_{n=1}^{N} [t_n - \boldsymbol{w}^\top \Phi_k(x_n)]^2}{2\sigma^2}$$

and the log of the posterior is

$$\ln p(\boldsymbol{w}|\mathcal{H}_k) = C - \frac{\boldsymbol{w}^\top \boldsymbol{w}}{2\sigma_w^2},$$

so that

$$A_{jl} = \frac{\sum_{n=1}^{N} [\Phi_k(x_n)]_j\, [\Phi_k(x_n)]_l}{\sigma^2} + \frac{\delta_{jl}}{\sigma_w^2}.$$

Hence:

$$A = \frac{\Phi^\top \Phi}{\sigma^2} + \frac{\mathrm{Id}}{\sigma_w^2}$$

where $[\Phi_k(x_n)]_l$ is the $l$-th entry of the vector $\Phi_k(x_n)$ and $\delta_{jl}$ is the Kronecker delta function. We see that in this case $A$ is independent of the choice $\boldsymbol{w} = \boldsymbol{w}_{\mathrm{MAP}}$. We can now use the $\boldsymbol{w}_{\mathrm{MAP}}$ obtained iteratively to compute the numerator of the evidence approximation.

Notice that in this particular case the structure of $A$ is such that we can actually do model selection iteratively as $A_{jl}$ can be easily obtained from a smaller dataset (for example $N-1$ instead of $N$).

## 7. Simulation results

We now discuss three examples of the results shown by our code.[3] To this purpose we always use two plots: the top one displays the whole data as red dots and the fitted model functions in lines of various colors. We can already see that some models fit better than others. However, to avoid overfitting we will not always use the best-fitting model but instead rely on Bayesian model selection.

The bottom plot shows the probability for each of the models given the data up to the time given on the *x*-axis. The red dots appear in exactly the order in which they go from left to right. This means that a specific point in the bottom plot corresponds to the whole range of intermediate data up to that position.

### 7.1. Low noise, quadratic growth data.

Figure 4 shows an example of a slight quadratic growth. The first 10 data points could be attributed to a constant model and some noise, so Bayesian model selection picks that (the blue curve is on top). After another 10 data points, the constant model cannot hold anymore so there is a fight between order 1 and order 2 polynomials. Order 1 polynomials are "simpler" but order 2 polynomials fit the slight convexity better, hence the algorithm is indifferent for some time until the convexity becomes too big (around datapoint no. 25) to ignore and model Poly2 wins finally.

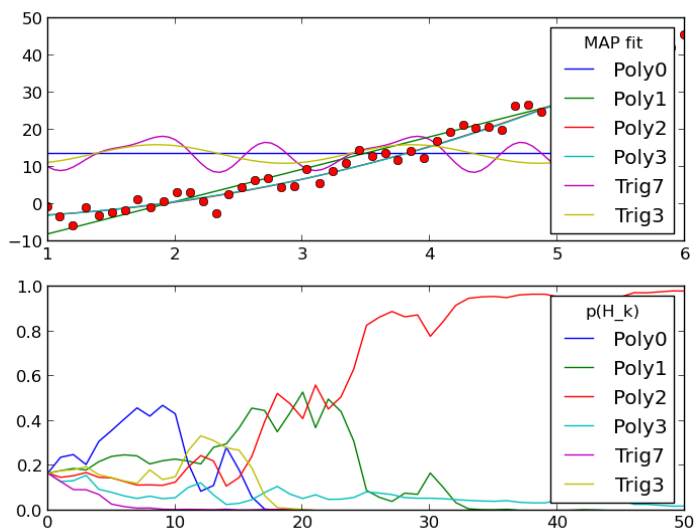

**Figure 4.** Data generated by a polynomial of order 2.

---

[3]. We implemented the algorithms discussed in Python using numpy and matplotlib. The code as well as the source for this document are freely available at github.com/mdbenito/ModelSelection.

**7.2. High noise, long range oscillating data.**

In the dataset of figure 5 we see a sinusoidal function generating the data (Trig3). The optical match is best for Trig3 and Trig7 although the latter is not preferred by the algorithm as it is too complicated and gets penalized by the Occam factor. In the bottom plot we see the algorithm choosing first a constant model which is very simple and explains the data well enough. As soon as the data shows consistent oscillation with some particular frequency (this is actually even hard to see with the naked eye due to the noise disturbing the result!), the algorithm correctly classifies the data as being generated by a lower order trigonometric function.

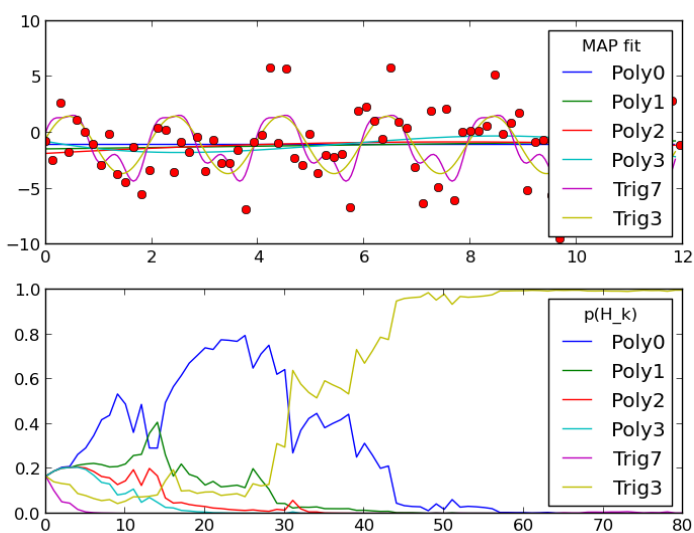

**Figure 5.** Trigonometric data generator, high noise and long time-scale.

**7.3. Low noise, short range oscillatory data.**

In figure 6 we see another Trig3 model generating data but with a lot less noise. This leads to another interesting phenomenon: After a short timespan where the constant model is deemed plausible, the linear polynomial Poly1 becomes popular as it explains the linear growth in the first quarter of the sinusoidal function's cycle well. As soon as the oscillation moves back down the quadratic Polynomial Poly2 is better suited for the data, but the next change of direction strengthens the Trig3 hypothesis. Due to a particular strange noise event around the 45th data point, Trig3, Poly2 and Poly3 are all seen as approximately similarly likely until the full cycle of the sinusoidal conclusively elects Trig3 as the most probable hypothesis.

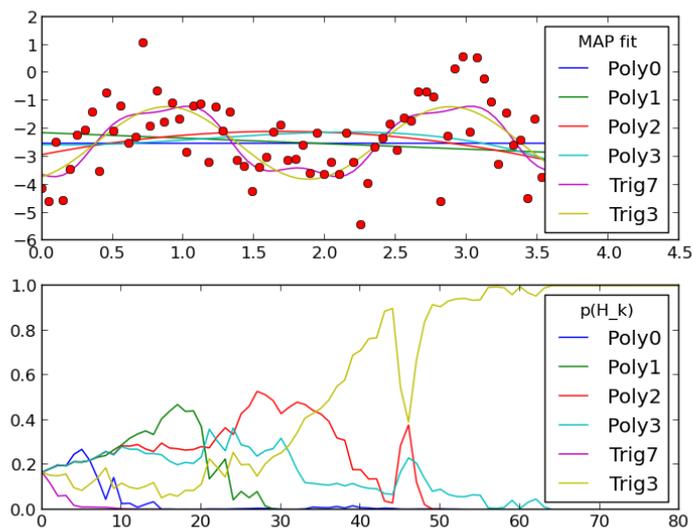

**Figure 6.** Trigonometric data generator, low noise and long time-scale

## Bibliography

**[Bis06]** Christopher M. Bishop. *Pattern recognition and machine learning*. Information science and statistics. Springer, 1 edition, aug 2006.

**[Mac05]** David J.C. MacKay. *Information theory, inference, and learning algorithms*. Cambridge University Press, Version 7.2 (4th printing) edition, mar 2005.